\documentclass[twoside,12pt,leqno]{article}
\advance\topmargin by -\headheight\advance\topmargin by -\headsep

\usepackage{amsmath}
\usepackage{amsthm}
\usepackage{amssymb}
\usepackage{amscd}
\usepackage{graphicx}
\usepackage{epsfig}
\numberwithin{equation}{section}
\newtheorem{theorem}{Theorem}[section]

\newtheorem{proposition}[theorem]{Proposition}

\theoremstyle{definition}
\newtheorem{definition}[theorem]{Definition}

\numberwithin{equation}{section}

\allowdisplaybreaks
\begin{document}
\thispagestyle{empty}




\bigskip

\bigskip

\centerline{{\Large  Bernoulli-Taylor formula in the case of
$Q$-umbral Calculus }}

\bigskip
\bigskip
\centerline{{\large by Ewa Krot-Sieniawska}}

\vspace*{.7cm}

\begin{abstract}

In this note we derive the $Q$-difference Bernoulli-Taylor formula
with the rest term of the Cauchy form by the Viskov's method
\cite{1}. This is an extension of technique presented in
\cite{5,6, 7} by the use of $Q$-extented Kwa\'sniewski's
$*_{\psi}$-product \cite{3, 4}. The main theorems of $Q$-umbral
calculus were given by G. Markowsky in 1968 (see \cite{2}) and
extented by A.K.Kwa\'sniewski \cite{3,4}.
\end{abstract}

KEY WORDS: Generalized differential operator, $*_{Q}$-product,
$Q$-integral

AMS 2000  numbers: 17B01, 17B35, 33C45,

\pagestyle{myheadings}
\section{Introduction - $Q$-umbral calculus}

We shall denote by {\bf P} the algebra of polynomials over the
field {\bf F}
 of characteristic zero.

 Let us consider a one parameter family \cal{F} of seqences. Then a sequence
  $\psi$ is called admissible (\cite{3,4}) if $\psi\in\cal{F}$. Where
\begin{multline*}
\mathcal{F}=\{\Psi:\mathbf{R}\supset[a,b];q\in[a,b]:\Psi(q):\mathbf{Z}
\rightarrow\mathbf{F};
\Psi_{0}(q)=1,\Psi_{n}(q)\neq0,\\
\Psi_{-n}(q)=0, n\in {\mathbf{N}}\}
\end{multline*}
Now let us to introduce the $\Psi$-notation \cite{3,4}:
$$ n_{\psi}=\Psi_{n-1}(q)\Psi_{n}^{-1}(q);$$
$$n_{\psi}!=n_{\psi}(n-1)_{\psi}\cdots2_{\psi}1_{\psi}=\Psi^{-1}_{n}(q),$$
$$n^{\underline{k}}_{\psi}= n_{\psi}(n-1)_{\psi}\cdots(n-k+1)_{\psi},$$
$$ \binom{n}{k}_{\psi}=\frac{ n^{\underline{k}}_{\psi}}{k_{\psi}!},$$
$$ exp_{\psi}\{y\}=\sum_{k=0}^{\infty}\frac{y^{k}}{k_{\psi}!}$$
\begin{definition}\cite{2,3}
Let $Q$ be a linear map $Q:{\bf P}\rightarrow {\bf P}$ such that:
$$\forall p\in {\bf{P}} \;\;\;\;deg(Qp)=(degp)-1$$
$degp=-1$ means $p=const=0$. Then $Q$ is called a generalized
differential operator \cite{2}.
\end{definition}

\begin{definition} \cite{2}
A normal sequence of polynomials $\{q_{n}\}_{n \geq 0}$ has the
following
 properties:
 \renewcommand{\labelenumi}{(\alph{enumi})}
\begin{enumerate}
\item $degq_{n}(x)=n$; \item $\forall x \in {\bf F}\; \; \; \;
q_{0}(x)=1$ ; \item $\forall n\geq 1 \; \; \; \; q_{n}(0)=0$ .
\end{enumerate}
\end{definition}

\begin{definition}\cite{2,3}
Let ${q_{n}(x)}_{n \geq 0}$ be the normal sequence. Then we call
it the $\psi$-basic sequence of the generalized differential
operator $Q$ (or
 $Q-\psi$-basic sequence) if:
 $$\forall n\geq 0\; \; \; \; Qq_{n}(x)=n_{\psi}q_{n-1}(x)$$
\end{definition}

In \cite{2} it is shown that once a differential operator $Q$ is
given a unique $\psi$-basic polynomial sequence and the other way
round: given a normal sequence ${q_{n}(x)}_{n\geq 0}$ there exists
a uniquely determined generalized differential operator $Q$.

\begin{definition} \cite{2,3}
The $\hat{x}_{Q}$-operator ($Q$-multiplication operator, the
operator dual
 to $Q$ ) is the linear map $\hat{x}_{Q}:{\bf{P}}\rightarrow {\bf{P}}$ such
 that:
 $$\forall n\geq 0 \; \; \; \; \hat{x}_{Q}q_{n}(x)=\frac{n+1}{(n+1)_{\psi}}q_{n+1}(x)$$
 Note that: \; \; $[Q,\hat{x}_{Q}]=id$
 \end{definition}
 \begin{definition} \cite{3}
 Let ${q_{n}(x)}_{n\geq 0}$ be a $Q-\psi$-basic sequence. Let
 $$E^{y}_{q}(Q)=E^{y}(Q)=exp_{Q,\psi}{yQ}=\sum_{k=0}^{\infty}
 \frac{q_{k}(y)Q^{k}}{k_{\psi}!}$$
 $E^{y}_{q}(Q)=E^{y}(Q)$ is called the $Q-\psi$-generalized translation
  operator.
\end{definition}

As was annouced in \cite{5,7}, the notion of Kwa\'sniewski's
$*_{\psi}$ product and its properties presented in \cite{3} can be
easily  $Q$-exteted as follows.
\begin{definition}\cite{3}
$$x*_{Q}q_{n}(x)=\hat{x}_{Q}(q_{n}(x))=\frac{n+1}{(n+1)_{\psi}}q_{n+1}(x), \
\ \ n\geq 0$$
$$x^{n}*_{Q}q_{n}(x)=(\hat{x}^{n}_{Q})(q_{1}(x))=
\frac{(n+1)!}{(n+1)_{\psi}!}q_{n+1}(x), \ \ \  n\geq 0$$ Therefore
$$x*_{Q}\alpha 1=x*_{Q}\alpha q_{0}(x)=\hat{x}_{Q}(\alpha q_{0}(x))=
\alpha \hat{x}_{Q}(q_{0}(x))=\alpha x*_{Q}1$$ and
$$f(x)*_{Q}q_{n}(x)=f(\hat{x}_{Q})q_{n}(x)$$
for every formal series f.
\end{definition}

\begin{definition}

According to definition above and \cite{3} we can define
$Q$-powers of $x$ by reccurence
 relation:
 $$x^{0*_{Q}}=1=q_{0}(x)$$
 $$x^{n*_{Q}}=x*_{Q}(x^{(n-1)*_{Q}})=\hat{x}_{Q}(x^{(n-1)*_{Q}})$$
 \end{definition}
 It is easy to show that:
 $$x^{n*_{Q}}=x*_{Q}x*_{Q}\ldots *_{Q}1=\frac{n!}{n_{\psi}!}q_{n}(x), \ \ \
 n\geq 0$$
 Also note that:
 $$x^{n*_{Q}}*_{Q}x^{k*_{Q}}=\frac{n!}{n_{\psi}!}q_{k+n}(x) $$
 and
 $$x^{k*_{Q}}*_{Q}x^{n*_{Q}}=\frac{k!}{k_{\psi}!}q_{k+n}(x) $$
 so in general i.e. for arbitrary admissible $\psi$ and for every
 $\{q_{n}\}_{n \geq 0}$ it is noncommutative.

 Due to above definition one can proof the following $Q$-extented
 properties of Kwa\'sniewski's $*_{\psi}$ product \cite{3,4}.
 \begin{proposition} \label{obstwoone} {\em   Let $f,g$ be formal series.
Then for $*_{Q}$ defined above holds:
\renewcommand{\labelenumi}{(\alph{enumi})}
\begin{enumerate}
\item $Qx^{n*_{Q}}=nx^{(n-1)*_{Q}}, \;\; n\geq0$; \item
$\exp_{Q,\psi}[\alpha x]=\exp{ \alpha\hat{x}_{Q}}1$   where
$\exp_{Q,\psi}{\alpha x}=\sum_{k\geq
0}\frac{q_{k}(x)\alpha^{k}}{k_{\psi}!};$ \item
$Q(x^{k}*_{Q}x^{n*_{Q}})=
(Dx^{k})*_{Q}x^{n*_{Q}}+x^{k}*_{Q}(Qx^{n*{Q}})$; \item
$Q(f*_{Q}g)=(Df)*_{Q}g+f*_{Q}(Qg),$ ($Q$-Leibnitz rule); \item
$f(\hat{x}_{Q})g(\hat{x}_{Q})1=f(x)*_{Q}\tilde{g};\;\;
\tilde{g}(x)=g(\hat{x}_{Q})1$.
\end{enumerate}}
\end{proposition}
According to \cite{2,3}, let us to define $Q$-integration operator
which is a right inverse operation to generalized differential
operator $Q$ i.e.:
$$Q\circ\int d_{Q}t=id$$
\begin{definition}
We define $Q$-integral as a linear operator such that
$$\int q_{n}(x)d_{Q}x=\frac{1}{(n+1)_{\psi}}q_{n+1}(x);\;\;\;n \geq 0.$$
\end{definition}
\begin{proposition} {\em
 \renewcommand{\labelenumi}{(\alph{enumi})}
\textrm{   }\\
\begin{enumerate}
\item $Q\circ \int_{\alpha}^{x}f(t)d_{Q}t=f(x);$ \item
$\int_{\alpha}^{x}(Qf)(t)d_{Q}t=f(x)-f(\alpha);$ \item formula for
integration "per partes" :
$$\int_{\alpha}^{\beta}(f*_{Q}Qg)(x)d_{Q}x=
[(f*_{Q}g)(x)]^{\beta}_{\alpha}-\int_{\alpha}^{\beta}((Df)*_{Q}g)(x)d_{Q}x.$$
\end{enumerate}  }
\end{proposition}

\section{$Q$-umbral calculus Bernoulli-Taylor formula}
In \cite{1} O.V.Viskov establishes the following identity
\begin{equation} \label{2.1}
\hat{p}\sum_{k=0}^{n} \frac{(-\hat{q})^{k}\hat{p}^{k}}{k!}=
\frac{(-\hat{q})^{n}\hat{p}^{n+1}}{n!}.
\end{equation}
what he calles the Bernoulli identity. Here $\hat{p},\hat{q}$
stand for linear operators satysfying condition:
\begin{equation} \label{2.2}
[\hat{p},\hat{q}]=\hat{p}\hat{q}-\hat{q}\hat{p}=id
\end{equation}
Now let $\hat{p}$ and $\hat{q}$ be as below:
$$\hat{p}=Q,\;\;\hat{q}=\hat{x}_{Q}-y,\;\; y\in \mathbf{F}$$
 From definition (1.4) we have $[Q,\hat{x}_{Q}-y]=id$.
 After submition into (\ref{2.1})  we get:
$$ Q\sum_{k=0}^{n}\frac{(y-\hat{x}_{Q})^{k}Q^{k}}{k!}=
\frac{(y-\hat{x}_{Q})^{n} Q^{n+1}}{n!}$$ Now let us apply it to
any polynomial (formal series) $f$:
$$ Q\sum_{k=0}^{n}\frac{(y-\hat{x}_{Q})^{k}(Q^{k}f)(t)}{k!}=
\frac{(y-\hat{x}_{Q})^{n} (Q^{n+1}f)(t)}{n!}.$$ Now using
definitions (1.6) and (1.7) of $*_{Q}$-product we get:
$$Q\sum_{k=0}^{n}\frac{(y-x)^{k*_{Q}}*_{Q}(Q
^{k}f)(t)}{k!}=\frac{(y-x)^{n*_{Q}}*_{Q}(Q^{n+1}f)(t)}{n!}.$$
After integration $\int_{y}^{x}d_{Q}t$ using proposition (1.2) it
gives $Q$-difference calculus Bernoulli-Taylor formula of the
form:
\begin{equation} \label{2.3}
f(x)=\sum_{k=0}^{n}\frac{1}{k!}(x-y)^{k*_{Q}}*_{Q}
(Q^{k}f)(y)+R_{n+1}(x)
\end{equation}
where $R_{n+1}$ stands for the rest term of the Cauchy type :
\begin{equation} \label{2.4}
R_{n+1}(x)=\frac{1}{n!}\int_{y}^{x}(x-t)^{n*_{Q}}*_{Q}
(Q^{n+1}f)(t)d_{\psi}t.
\end{equation}
\section{Special cases}
\begin{enumerate}
\item An example of generalized differential operator is $Q\equiv
 \partial_{\psi}$. Then $q_{n}=x^{n}$ and $\partial_{\psi}x^{n}=
 n_{\psi}x^{n-1}$ for an admissible $\psi$. $\partial_{\psi}$ is called
 $\psi$-derivative. Then also $\hat{x}_{Q}\equiv \hat{x}_{\psi}$ and
 $\hat{x}_{\psi}x^{n}=\frac{(n+1)}{(n+1)_{\psi}}x^{n+1}$ and
 $[\partial_{\psi}, \hat{x}_{\psi}]=id.$ In this case we get $\partial_{\psi}$
 -difference calculus Bernoulli-Taylor formula presented in \cite{5,6,7} of
  the form
$$ f(x)=\sum_{k=0}^{n}\frac{1}{k!}(x-\alpha)^{k*_{\psi}}*_{\psi}
(\partial_{\psi}^{k}f)(\alpha)+R_{n+1}(x)$$ with
$$R_{n+1}(x)=\frac{1}{n!}\int_{\alpha}^{x}(x-t)^{n*_{\psi}}*_{\psi}
(\partial_{\psi}^{n+1}f)(t)d_{\psi}t.$$ In \cite{3} there is given
a condition for the case $Q=Q(\partial_{\psi})$
for some admissible $\psi$(see Section 2, Observation 2.1). \\
\item For $Q=\partial_{\psi}, \ q_{n}(x)=x^{n}, \ n\geq 0$ the
choice $\psi_{n}(q)=\frac{1}{[R(q^{n})]!}, \ R(x)=\frac{1-x}{1-q}$
gives the well known $q$-factorial $n_{q}!=n_{q}(n-1)_{q}\ldots
2_{q}1_{q}$, for $n_{q}=1+q+q^{2}+\ldots +q^{n-1}$. Then
$\partial_{\psi}=\partial_{q}$ becomes the well known Jackson's
derivative $\partial_{q}$:
$$(\partial_{q}f)(x)=\frac{f(x)-f(qx)}{(1-q)x}$$
The $\partial_{q}$-difference version of the Bernoulli-Taylor
formula was
given in \cite{5} by the use of $*_{q}$-product.\\
\item By the choice $Q\equiv D\equiv \frac{d}{dx}, \
q_{n}(x)=x^{n}$ and
 $\psi_{n}=\frac{1}{n!}$ after submition to (\ref{2.3}), (\ref{2.4}) we get
  the classical Bernoulli-Taylor formula of the form:
  $$f(x)=\sum_{k=0}^{n}\frac{(x-\alpha)^{k}}{k!}f^{(k)}(\alpha)+\int_{\alpha}
  ^{x}\frac{(x-t)^{n}}{n!}f^{n+1}(t)dt$$
  where $f^{(k)}(\alpha)=(D^{k}f)(\alpha)$. \\
\end{enumerate}
{\bf {Acknowledgements}}

I would like to thank Professor A. Krzysztof Kwa´sniewski for his
very helpful comments, suggestions, improvements and corrections
of this note.


\noindent
      Ewa Krot-Sieniawska \\
      Institute of Computer Science,  University in Bia{\l}ystok\\
      PL-15-887 Bia{\l}ystok, ul.Sosnowa 64, POLAND\\
      e-mail: ewakrot@wp.pl

\label{lastpage}
\end{document}